\documentclass{segabs}

\usepackage{amsmath,amsfonts,amssymb,mathrsfs}
\usepackage{xcolor}

\usepackage{graphicx}
\graphicspath{{./FIGURES/}}

\newtheorem{thm}{Theorem}
\newtheorem{algorithm}[thm]{Algorithm}

\DeclareMathAlphabet{\itbf}{OML}{cmm}{b}{it}
\DeclareMathAlphabet\mathbfcal{OMS}{cmsy}{b}{n}

\renewcommand{\hat}{\widehat}

\def\RR{\mathbb{R}}
\def\bx{{{\itbf x}}}

\def\bu{{{\itbf u}}}

\def\br{{\itbf r}}

\def\bd{{\itbf d}}

\def\bU{{\itbf U}}
\def\bJ{{\itbf J}}

\def\bet{{\boldsymbol{\eta}}}
\def\ss{{(s)}}

\def\bV{{\itbf V}}
\def\bR{{\itbf R}}

\def\bI{{\itbf I}}
\def\bS{{\itbf S}}

\def\bM{{\itbf M}}
\def\bD{{\itbf D}}
\def\bA{\boldsymbol{\cal A}}

\def\FWI{{\scalebox{0.5}[0.4]{FWI}}}

\def\RM{{\scalebox{0.5}[0.4]{ROM}}}

\def\om{\omega}

\def\12{{\frac{1}{2}}}

\begin{document}

\title{Velocity estimation via model order reduction}

\renewcommand{\thefootnote}{\fnsymbol{footnote}} 

\author{Alexander V. Mamonov\footnotemark[1], University of Houston;
Liliana Borcea, University of Michigan; Josselin Garnier, Ecole Polytechnique;
and J\"{o}rn Zimmerling, University of Michigan}

\footer{Velocity estimation via model reduction}
\lefthead{A.V. Mamonov, L. Borcea, J. Garnier \& J. Zimmerling}
\righthead{Velocity estimation via model reduction}

\maketitle

\begin{abstract}
A novel approach to full waveform inversion (FWI), based on a data driven reduced order model (ROM) of the wave 
equation operator is introduced. The unknown medium is probed with pulses and the time domain pressure 
waveform data is recorded on an active array of sensors. The ROM, a projection of the wave equation operator 
is constructed from the data via a nonlinear process and is used for efficient velocity estimation. 
While the conventional FWI via nonlinear least-squares data fitting is challenging without 
low frequency information, and prone to getting stuck in local minima (cycle skipping), minimization of ROM 
misfit is behaved much better, even for a poor initial guess. For low-dimensional parametrizations of the 
unknown velocity the ROM misfit function is close to convex. The proposed approach consistently outperforms
conventional FWI in standard synthetic tests.
\end{abstract}

\section{Introduction}

We consider the inverse problem of velocity estimation from time-domain reflection data recorded by an 
array of $m$ sensors that can both emit and record. For simplicity we work with the acoustic wave equation
with unknown velocity $c(\bx)$. The proposed approach can be extended to vectorial (elastic) waves.

The model pressure wave $p^\ss(t,\bx)$ generated by the $s^{\rm th}$ sensor, for $s = 1, \ldots, m$, 
satisfies the initial value problem
\begin{align}
\big[ \partial_t^2 - c^2(\bx) \Delta \big]  p^\ss(t,\bx) &= f'(t) \theta(\bx-\bx_s), \quad 
t \in \RR, 
\label{eq:I1} \\
p^\ss(t,\bx) & = 0, \quad t <-t_f ,\label{eq:I2}
\end{align}
for $\bx \in \Omega$, a simply connected domain in two or three dimensions, with boundary $\partial \Omega$. 
We set homogeneous boundary conditions (Dirichlet, Neumann, or a combination thereof) on $\partial \Omega$.
We assume the sensors in the array to be identical, modeled by a function $\theta(\bx)$, with small support around 
the origin. 
Each sensor emits the probing pulse $f(t)$, supported on $(-t_f,t_f)$. 
For simplicity we take $f(t)$ to be an even function, with a non-negative Fourier transform
$ \hat f(\om) \ge 0 $ that is analytic.

The inverse problem is to estimate the velocity $c(\bx)$ from the measurements
\begin{equation}
\mathcal{M}^{(r,s)}(t) = \int_\Omega d \bx \, \theta(\bx-\bx_r)p^\ss(t,\bx) ,\label{eq:I3}
\end{equation}
for $s,r = 1,\ldots, m$ and $t \in [0,T]$. 
Conventional FWI approach to velocity estimation, see, e.g., 
\cite{tarantola1984inversion}, 
suffers from a fundamental flaw: the objective function is nonconvex with numerous local minima.  
This effect makes any local optimization algorithm unlikely to succeed, 
in the absence of an accurate starting guess, see \cite{santosa1989analysis}. 
We address this issue by reformulating the optimization problem for velocity estimation using ROM.

\section{Theory}

\subsection{Symmetrized wave operator and data model}

It is convenient to work with a self-adjoint wave operator
\begin{equation}
{\cal A} = -c(\bx) \Delta \big[ c(\bx) \cdot \big],
\label{eq:defA}
\end{equation} 
a similarity transformation of $-c^2(\bx) \Delta$, the wave operator of equation~\ref{eq:I1}. 
To obtain the ROM of ${\cal A}$, we introduce the even in time wave
\begin{equation}
{\rm w}^\ss(t,\bx) =  \frac{\big[p^\ss(t,\bx) + p^\ss(-t,\bx)]}{c(\bx)},
\label{eq:M2}
\end{equation}
and ``the data", an $m \times m$ matrix $\bD(t)$, with entries 
\begin{equation}
D^{(r,s)}(t) = \frac{\mathcal{M}^{(r,s)}(t)+\mathcal{M}^{(r,s)}(-t)}{c(\bx_r)c(\bx_s)}, 
\label{eq:M3Dat}
\end{equation} 
for $s,r = 1, \ldots, m$, that can be obtained from the measurements $\mathcal{M}^{(r,s)}(t)$, 
assuming $c(\bx)$ is known in the vicinity of sensor center points $\bx_s = 1,\ldots,m$.

The ROM is computed from $2n-1$ equidistant time samples of the data matrix $\bD(t)$ and its second derivative
\begin{equation}
\bD_j ~ \mbox{and} ~~\ddot{\bD}_j = \partial_t^2 \bD(j \tau), \quad j = 0, \ldots, 2n-2,
\label{eq:M5}
\end{equation}
where the second derivative can be obtained from $\bD(t)$ via Fourier domain differentiation.
The sampling interval $\tau$ should be chosen according to the Nyquist sampling rate for 
the essential frequency of $f(t)$ (the largest frequency in the interval outside of which $\hat f(\om)$ is small).

Using the above, we rewrite Equations~\ref{eq:I1}--\ref{eq:I2} as the initial value problem
\begin{align}
&\big[\partial_t^2 + {\cal A} \big] u^\ss(t,\bx) = 0, \quad t > 0, ~\bx \in \Omega, \label{eq:M8} \\
&u^\ss(0,\bx)= u_0^\ss(\bx), ~~ \partial_t u^\ss(0,\bx) = 0, \quad \bx \in \Omega, \label{eq:M9}
\end{align}
where the wavefield
\begin{equation}
u^\ss(t,\bx) = \cos \big(t \sqrt{{\cal A}}\big) u_0^\ss(\bx), 
\label{eq:M7}
\end{equation}
is related to the even in time wave via 
\begin{equation}
{\rm w}^\ss(t,\bx) = \hat f^\12 \big(\sqrt{{\cal A}}\big) u^\ss(t,\bx),
\end{equation}
and the initial state is given by 
\begin{equation}
u_0^\ss(\bx) = \hat f^\12 \big(\sqrt{{\cal A}}\big) \frac{\theta(\bx-\bx_s)}{c(\bx_s)},
\quad s = 1, \ldots, m,
\label{eq:M10}
\end{equation}
which is supported in a ball centered at $\bx_s$, with radius of order $c(\bx_s)t_f$.
The details of the above formulation can be found in \cite{borcea2020reduced}.

Working with wavefields $u^\ss(t,\bx)$ allows to express the data samples in a symmetric inner product form
\begin{align}
D_j^{(r,s)} & = \int_{\Omega} d \bx \, u_0^{(r)}(\bx) 
\cos \big( j \tau \sqrt{{\cal A}} \big) u_0^\ss(\bx), \label{eq:M11} \\
\hspace{-0.07in} \ddot{D}_j^{(r,s)} & = - \int_{\Omega} d \bx \, u_0^{(r)}(\bx) {\cal A} 
\cos \big( j \tau \sqrt{{\cal A}} \big) u_0^\ss(\bx) , \label{eq:M12}
\end{align}
for $r,s = 1, \ldots, m$ and $j = 0, \ldots, 2n-2$. These formulas can be simplified using block algebra notation.  
Gathering all the waves $u^\ss(j \tau,\bx)$ into the row vector field called a snapshot
$ \bu_j(\bx) = \big(u^{(1)}(j\tau,\bx), \ldots, u^{(m)}(j\tau,\bx) \big) $, we observe that it satisfies 
\begin{equation}
\bu_j(\bx) = \cos \big(j \tau\sqrt{{\cal A}}\big) \bu_0(\bx), \quad j \ge 0,
\label{eq:M15}
\end{equation}
hence the data matrix samples can be written as
\begin{align}
\bD_j &= \langle \bu_0, \bu_j \rangle = \langle \bu_0,  \cos \big(j \tau  \sqrt{{\cal A}}\big) \bu_0 \rangle,
\label{eq:M16} \\
\ddot \bD_j &=  - \langle \bu_0, {\cal A} \bu_j \rangle,
\label{eq:M17}
\end{align}
where we denote by
$\langle \boldsymbol{\phi},\boldsymbol{\psi} \rangle = \int_\Omega d \bx \, 
\boldsymbol{\phi}^T(\bx) \boldsymbol{\psi}(\bx)$
the integral of the outer product of any functions $\boldsymbol{\phi}(\bx)$ and 
$\boldsymbol{\psi}(\bx)$ with values in $\RR^{1\times m}$ or $\RR^{1\times nm}$ 
and $T$ stands for the transpose. 

\subsection{Wave operator ROM}

At the core of the proposed approach is the ROM of the wave operator, the orthogonal projection
of ${\cal A}$ onto the space
\begin{equation}
\mathbb{S} = \mbox{span} \big\{ \bu_0(\bx),\ldots, \bu_{n-1}(\bx) \big\},
\label{eq:M19}
\end{equation}
spanned by the first $n$ snapshots. Gathering these snapshots into the $nm$-dimensional row vector field 
\begin{equation}
\bU(\bx) = \big( \bu_0(\bx), \ldots, \bu_{n-1}(\bx) \big) \in \RR^{1 \times nm},
\label{eq:M23}
\end{equation}
an orthonormal basis for $\mathbb{S}$ can be obtained via the block Gram-Schmidt orthogonalization 
$\bU(\bx) = \bV(\bx) \bR$,
where $\left\langle \bV, \bV \right\rangle = \bI_{nm}$ and $\bR$ is an $nm \times nm$
block upper triangular matrix, a block Cholesky factor of the so-called mass matrix
\begin{equation}
\bM = \left\langle \bU, \bU \right\rangle = \bR^T\bR \in \RR^{nm \times nm}.
\end{equation}

A remarkable property of the mass matrix $\bM$ is that it can be computed from data samples only. 
Thus, it is possible to obtain the block Cholesky factor $\bR$ without the knowledge of internal wavefields
that is otherwise required for Gram-Schmidt orthogonalization. 
Explicitly, using the trigonometric identity
\begin{equation}
\label{eq:trig}
\cos(\alpha) \cos(\beta) = \frac{1}{2} \left[\cos(\alpha + \beta) + \cos(\alpha-\beta)\right],
\end{equation}
the $m \times m$ blocks of the mass matrix can be computed as
\begin{equation}
\bM_{i,j} = \langle \bu_i, \bu_j \rangle 
= \frac{1}{2} \left( \bD_{i+j} + \bD_{|i-j|} \right) \in \RR^{m \times m}, 
\label{eq:M32}
\end{equation}
for $i,j = 0,\ldots,n-1$. This computation is the first crucial step in computing the ROM of ${\cal A}$,
given by 
\begin{equation}
 \bA^\RM =  \langle \bV, \bA \bV \rangle  = 
 - \int_{\Omega} d\bx \, \bV^T(\bx) c(\bx) \Delta \big[ c(\bx) \bV(\bx) \big].
\label{eq:M39V}
\end{equation}
Indeed, substituting $\bV(\bx) = \bU(\bx) \bR^{-1}$ into the above expression, we obtain
\begin{equation}
\bA^\RM = \bR^{-T} \langle \bU, \bA \bU \rangle \bR^{-1}.
\end{equation}

Remarkably, the so-called $nm \times nm$ operator stiffness matrix $\bS = \langle \bU, \bA \bU \rangle$
can also be computed from the data samples using a calculation similar to equation~\ref{eq:M32} for the
mass matrix. Explicitly, the $m \times m$ blocks of $\bS$ are given by
\begin{equation}
\bS_{i,j} = \langle \bu_i, {\cal A} \bu_j \rangle 
= -\frac{1}{2} \left(\ddot \bD_{i+j} + \ddot \bD_{|i-j|} \right) \in \RR^{m \times m},
\label{eq:M43}
\end{equation}
for $i,j = 0, \ldots, n-1$. This explains the need for the second derivative data samples. 
We summarize the computation of the wave operator ROM $\bA^\RM$ in the following algorithm.

\begin{algorithm}\textbf{\emph{(Data-driven ROM computation)}}
\label{alg:arom}
~\\
\noindent \textbf{Input:} The measurements $\mathcal{M}^{(r,s)}(t)$ for $t \in [0, T]$.\\
1. Compute $\{\bD_j, \ddot \bD_j\}_{j=0}^{2n-2}$ using equation~\ref{eq:M3Dat} and Fourier domain differentiation. \\
2. Calculate blocks of mass and stiffness matrices $\bM, \bS \in \mathbb{R}^{nm \times nm}$ 
using Equations~\ref{eq:M32} and \ref{eq:M43}, respectively. \\
3.  Perform the block Cholesky factorization $\bM = \bR^T \bR$ using, e.g., Algorithm 5.2 from
\cite{druskin2018nonlinear}. \\
\textbf{Output:} wave operator ROM $\bA^\RM = \bR^{-T} \bS \bR^{-1} $.
\end{algorithm}

\subsection{ROM based velocity estimation}

The proposed approach for velocity estimation is based on minimizing the misfit of operator ROM instead 
of data misfit. We expect this to outperform the conventional FWI approach due to
a simple dependency of $\bA^\RM$ on the velocity. In particular, for a fixed projection space $\mathbb{S}$, 
and thus a fixed basis $\bV(\bx)$, the dependency of $\bA^\RM$ on $c(\bx)$ is quadratic, as observed in 
equation~\ref{eq:M39V}. While $\bV(\bx)$ also depends on $c(\bx)$, the numerical examples presented
below demonstrate that ROM misfit formulation leads to optimization objective that is 
close to convex. 

Given the above, the basic ROM based velocity estimation is to minimize the misfit of the wave operator ROM
\begin{equation}
\min_{v \in \mathcal{C}} \mathcal{O}^\RM (v), \quad 
\mathcal{O}^\RM (v) = \left\| {\rm Triu} \big( \bA^\RM(v) -  \bA^\RM \big) \right\|_2^2,
\label{eq:M46}
\end{equation}
where ${\rm Triu}$ extracts the upper triangular part (including the main diagonal) of a symmetric matrix
and stacks the entries in a vector. Hereafter $\| \cdot \|_2$ is the Euclidean norm.
Here $v$ denotes a velocity model in the search space $\mathcal{C}$ parametrized using appropriately 
chosen basis functions $\{\phi_l(\bx)\}_{l=1}^N$:
\begin{equation}
{v}(\bx;\bet) = c_o(\bx) + \sum_{l = 1}^N \eta_l \phi_l(\bx),
\label{eq:IM2}
\end{equation}
where $c_o(\bx)$ is the initial guess. Then, the optimization is for the vector 
$\bet = (\eta_1, \ldots, \eta_N)^T \in \RR^N$. The ROM $\bA^\RM$ is computed from the measurements 
with Algorithm~\ref{alg:arom}, and $\bA^\RM(v)$ is computed with the same algorithm for the measurements
calculated for the velocity model $v$. 

In practice, ROM based velocity estimation benefits from a layer stripping formulation. Then, instead of 
working with the whole matrix $\bA^\RM$, we consider the restriction $\big[\bA^\RM\big]_k$, an upper 
left $km \times km$ submatrix of $\bA^\RM$, where $k \le n$ increases gradually as the velocity 
estimation progresses.
Due to causality of wave operator ROM, $\big[\bA^\RM\big]_k$ is only affected by the first $2k-1$ 
data samples, see \cite{borcea2022waveform} for details. Another modification to the basic formulation 
in equation~\ref{eq:M46} is to discard some of the upper triangular entries of $\bA^\RM$ from the misfit 
calculation to decrease the computational burden and storage requirements. We include only the first few 
$dm$ diagonals in the objective function calculation, where $d$ is an integer between $1$ and $k$. 
We denote by ${\rm Rest}_{d,k}: \RR^{km \times km} \mapsto \RR^{dm (km- (dm-1)/2)}$
the mapping that takes a $km \times km $ matrix,  keeps only its first $dm$ upper diagonals, 
including the main one, and puts their entries into a column vector. We use the resulting objective function
\begin{equation}
\mathcal{O}_{d,k}({v}) = \left \|\mbox{Rest}_{d,k} \big(\left[\bA^\RM(v)-\bA^\RM\right]_k\big)\right\|_2^2,
\label{eq:newObj}
\end{equation}
to formulate the following velocity estimation algorithm.

\begin{algorithm}\textbf{\emph{(ROM based velocity estimation)}}
\label{alg:prowi}
~\\
\textbf{Input:} The ROM $\bA^\RM$ computed from the measurements. \\
1. Set the number $L$ of layers for the layer stripping approach and the number $q$ of iterations per layer.\\
2. Choose $L$ positive integers $\{k_l\}_{l=1}^L$, satisfying
\[1 \le k_1 \leq k_2 \leq  \cdots \leq  k_{L} = n.\] 
3. Starting with the initial vector $\bet^{(0)}= {\bf 0}$, proceed:\\
For $l = 1,2,\ldots,L$, and $j = 1,\ldots, q$, set the update index $i = (l-1)q+j$. 
Compute $\bet^{(i)}$ as a Gauss-Newton update for minimizing the functional
\begin{equation}
\mathcal{F}_i(\bet) = \mathcal{O}_{d,k_l} \big(v(\cdot; \bet)) + \mathcal{R}_i(\bet),
\end{equation}
linearized about $\bet^{(i-1)}$, where $\mathcal{R}_i(\bet)$ is a regularization penalty functional.\\
\textbf{Output:} velocity estimate $c^{\rm est}(\bx) = {v}(\bx;\bet^{(L q)})$.
\end{algorithm}


Even in the absence of noise, Algorithm~\ref{alg:prowi} requires the use of regularization term 
$\mathcal{R}_i(\bet)$. Note that in the presence of noise a more sophisticated regularization strategy is
needed, described in detail in \cite{borcea2022waveform}. We choose Tikhonov regularization 
$\mathcal{R}_i(\bet) = \mu_i \|\bet\|_2^2$,
where $\mu_i$ is chosen adaptively with the following procedure. Denote by
\begin{equation}
\mathcal{G}(\bet;d,k_l) = 
\mbox{Rest}_{d,k_l} \left(\left[\bA^\RM({v}(\cdot;\bet))-\bA^\RM\right]_{k_l}\right)
\label{eq:residual}
\end{equation}
the residual vector of objective $\mathcal{O}_{d,k_l} \big(v(\cdot; \bet))$.
Evaluated at $\bet = \bet^{(i-1)}$, its Jacobian is the matrix
\begin{equation}
\bJ^{(i)} = \nabla_\bet  \mathcal{G} (\bet^{(i-1)};d,k_l) \in \RR^{dm (km- (dm-1)/2) \times N}.
\end{equation}
We choose $N$ so that the Jacobian has more rows than columns.
If $\sigma_1^{(i)} \geq \sigma_2^{(i)} \geq \cdots \geq \sigma_{N}^{(i)}$
are the singular values of $\bJ^{(i)}$, then for a fixed parameter $\gamma \in (0.2, 0.4)$ 
(typically employed values, with smaller $\gamma$ corresponding to more regularization), we take
$ \mu_i = \big( \sigma^{(i)}_{\lfloor \gamma N \rfloor} \big)^2 $.

The regularized Gauss-Newton update direction is
\begin{equation}
\bd^{(i)} = - \left( \big( \bJ^{(i)} \big)^T \bJ^{(i)} + \mu_i \bI_N \right)^{-1}
\big( \bJ^{(i)} \big)^T \br^{(i)},
\end{equation}
where $\bI_N$ is the $N \times N$ identity matrix and $\br^{(i)} = \mathcal{G}(\bet^{(i-1)}; d, k_l)$.
Once the update direction $\bd^{(i)}$ is computed, we use a line search
\begin{equation}
\alpha^{(i)} = \mathop{\mbox{argmin}}\limits_{\alpha \in (0, \alpha_{\max})}
\mathcal{F}_i\big( \bet^{(i-1)} + \alpha \bd^{(i)} \big)
\end{equation}
to find the step length $\alpha^{(i)}$, where we take $\alpha_{\max} = 3$.
This gives the Gauss-Newton update $ \bet^{(i)} = \bet^{(i-1)} + \alpha^{(i)} \bd^{(i)} $.

\section{Examples}

We illustrate the performance of Algorithm~\ref{alg:prowi} and compare it to conventional FWI,
the minimization of data misfit
\begin{equation}
\min_{v \in \mathcal{C}} \mathcal{O}^\FWI({v}), \quad 
\mathcal{O}^\FWI({v}) = \sum_{j=0}^{2n-1} \left\| {\rm Triu} \big(\bD_j({v})-\bD_j \big)\right\|_2^2,
\label{eq:FWIObj}
\end{equation}
where $\bD_j({v})$ are the data matrix samples for model velocity.
All the results are for a band-limited source pulse
\begin{equation}
f(t) = \cos(\om_o t) \exp \Big[-\frac{(2 \pi B)^2 t^2}{2}\Big],
\label{eq:pulse}
\end{equation}
with central frequency  $\om_o/(2 \pi) = 6$Hz and bandwidth $B=4$Hz. 
The essential frequency is $\om_o/(2 \pi) + B = 10$Hz.

%

\subsection{Topography of objective functions}

\begin{figure}[ht!]
\begin{center}
\includegraphics[width=0.3\textwidth]{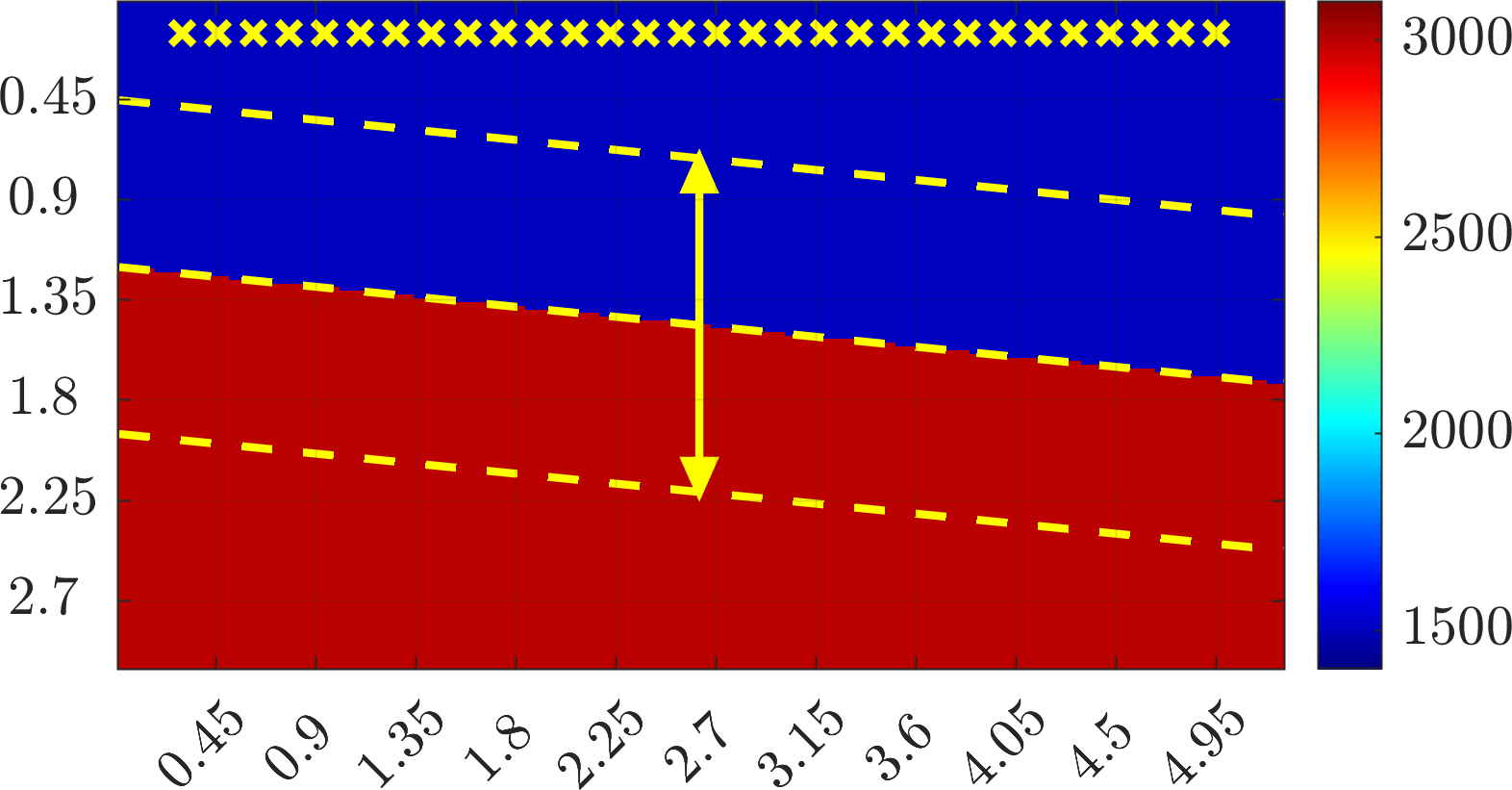}
\begin{tabular}{cc}
{\scriptsize Log of FWI objective $\mathcal{O}^\FWI$} & 
{\scriptsize Log of ROM objective $\mathcal{O}^\RM$} \\
\includegraphics[width=0.46\columnwidth]{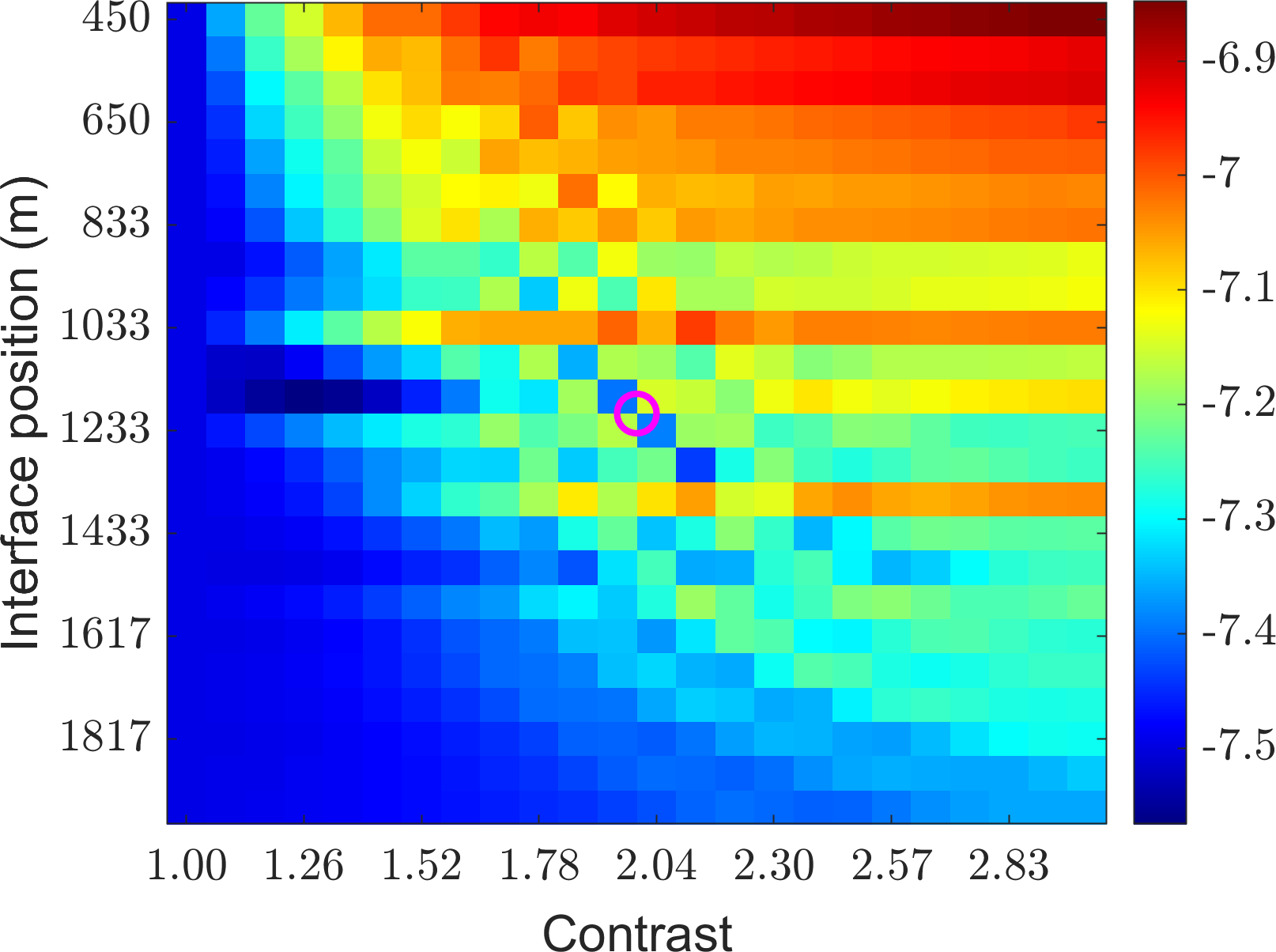} &
\includegraphics[width=0.46\columnwidth]{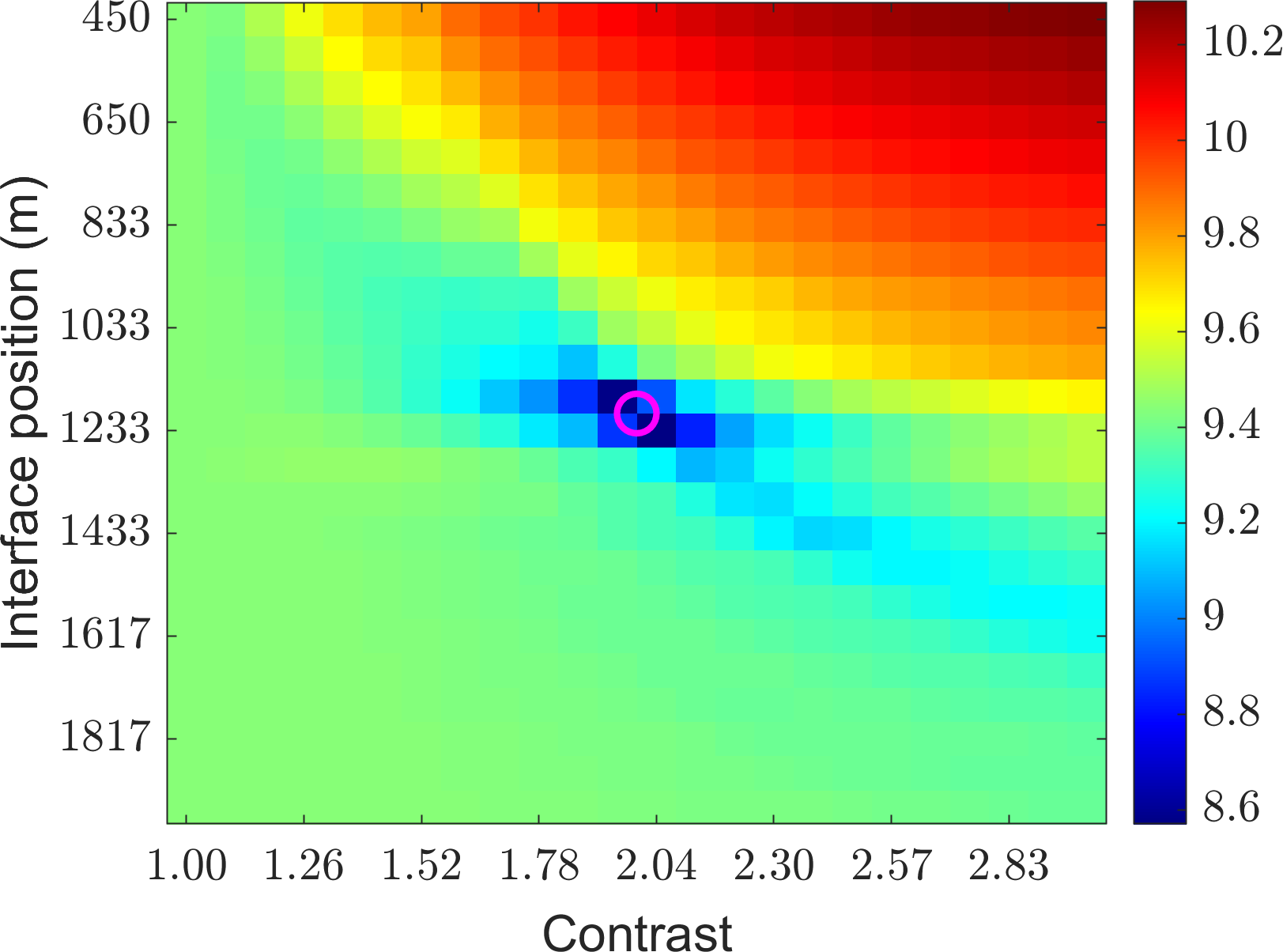}
\end{tabular}
\end{center}
\caption{Top: velocity model for objective topography study. Middle dashed line
shows the actual interface location; top and bottom dashed lines show
the extent of interface location parameter sweep. All $m=30$ sensors are
shown as yellow $\times$. Distances are in $\rm{km}$, velocity in $\rm{m/s}$.
Bottom: logarithms of objective functions.
True parameters are indicated with $\textcolor{magenta}{\bigcirc}$.}
\label{fig:topo}
\end{figure}

Consider the velocity model displayed in Figure~\ref{fig:topo}
consisting of two regions separated by a slanted interface
with velocities $c_t = 1500$m/s and $c_b = 3000$m/s above and below the interface, respectively.
To compare the objective functions, we sweep a two-parameter search space: 
the first parameter is the depth of the interface at its leftmost point
(actual depth is $1.2\rm{km}$); the second parameter is the contrast $c_b/c_t$
(actual contrast is $2$). 

In Figure~\ref{fig:topo} we compare the conventional FWI objective $\mathcal{O}^\FWI$ to
wave operator ROM objective $\mathcal{O}^\RM$. We observe that $\mathcal{O}^\FWI$ displays
numerous local minima, at points in the search space that are far from the true one. 
The horizontal stripes in FWI objective plot are manifestations of cycle skipping. 
The wave operator ROM objective is smooth and has a single minimum, at the true depth and contrast.
This confirms that the wave operator ROM misfit minimization is superior to conventional FWI 
formulation for velocity estimation since $\mathcal{O}^\RM$ is more friendly towards local
optimization algorithms.

\subsection{The ``Camembert" example}

\begin{figure}[ht!]
\begin{center}
\includegraphics[width=0.45\columnwidth]
{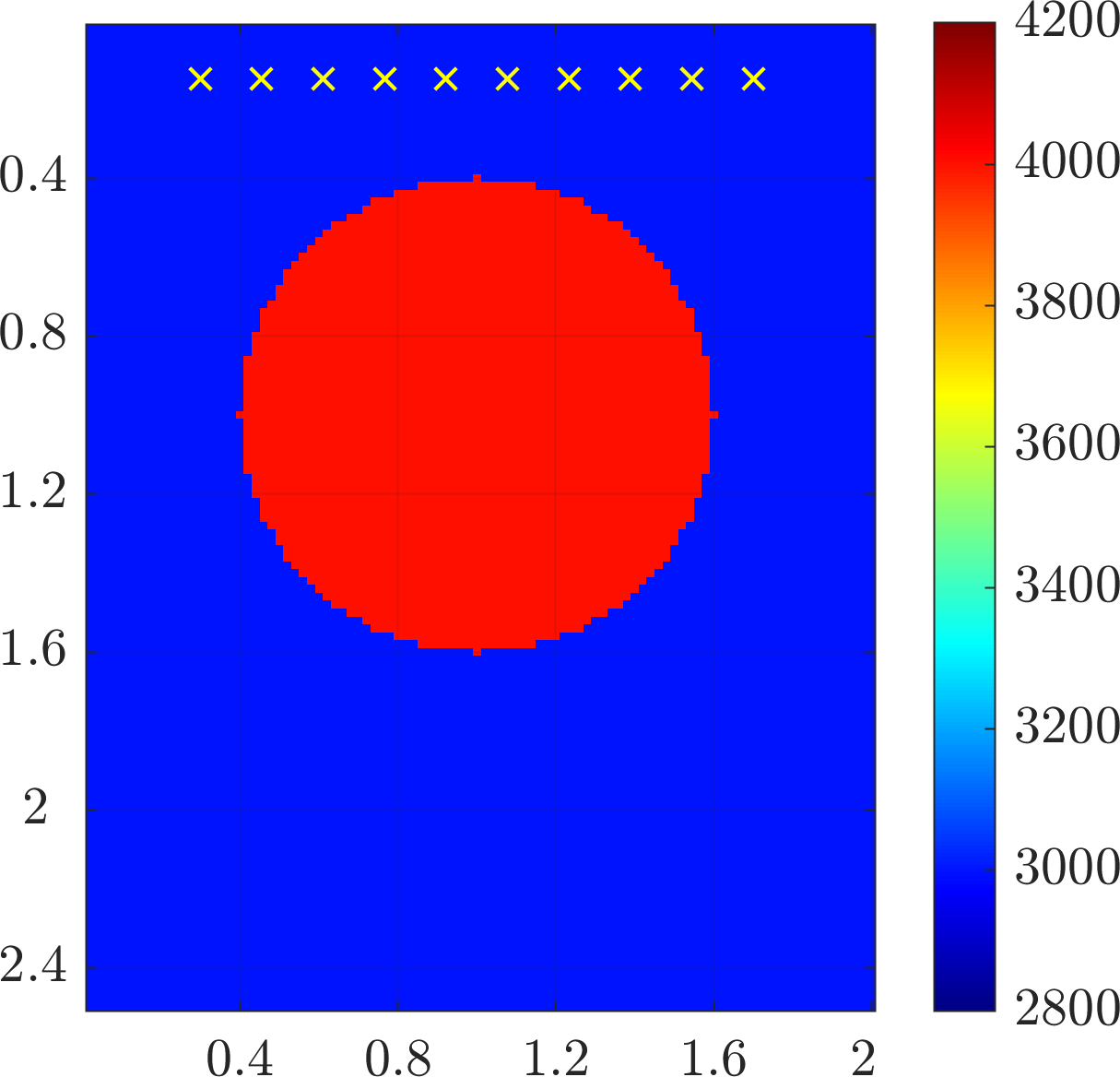}
\begin{tabular}{cc}
{\scriptsize Operator ROM estimate} &  {\scriptsize Conventional FWI estimate} \\
\includegraphics[width=0.45\columnwidth]
{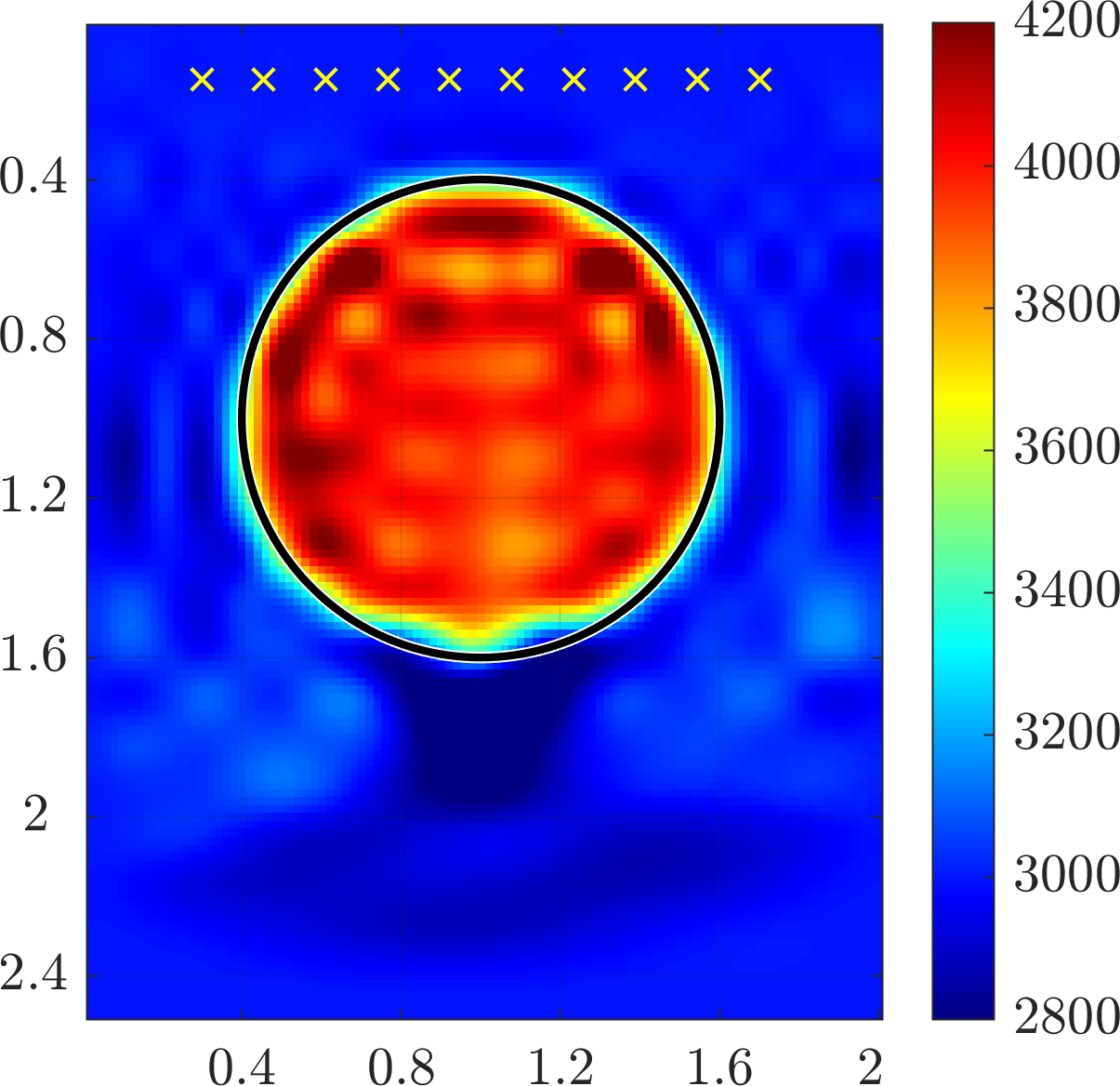} &
\includegraphics[width=0.45\columnwidth]
{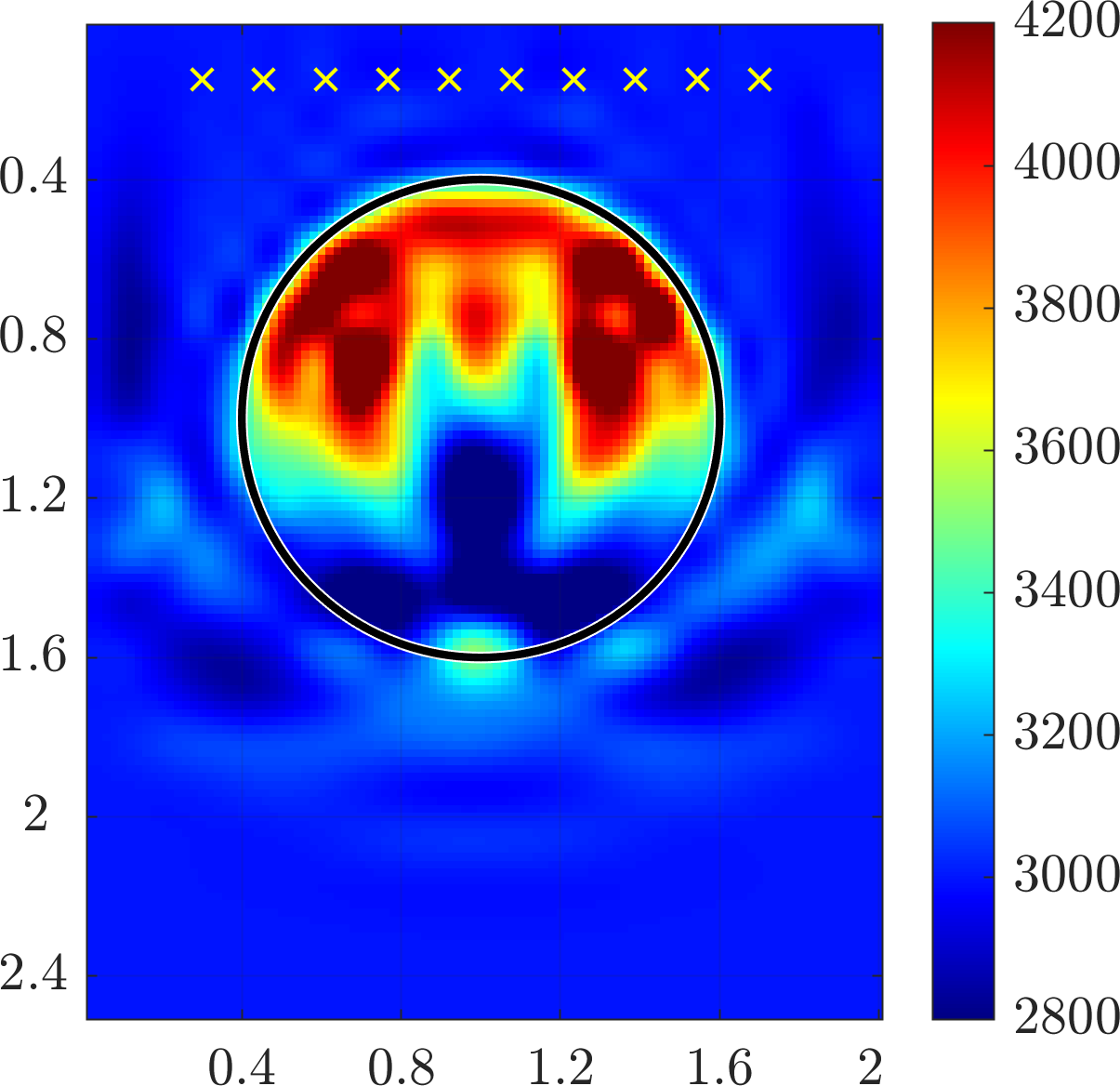}
\end{tabular}
\end{center}
\caption{Top: ``Camembert" velocity model. All $m=10$ sensors are shown as yellow $\times$. 
Bottom: estimated velocity after $60$ Gauss-Newton iterations, 
obtained with Algorithm~\ref{alg:prowi} (left) and conventional FWI (right). 
True inclusion boundary is shown as a black circle. 
The axes are in $\rm{km}$, velocity in $\rm{m/s}$.}
\label{fig:Camembert}
\end{figure}

Following \cite{yang2018application} we consider the ``Camembert" model with a circular inclusion 
of radius of $600$m, centered at $(1\rm{km},1\rm{km})$ in
$\Omega = [0, 2\rm{km}] \times [0, 2.5\rm{km}]$,
with $c(\bx) = 4000\rm{m/s}$ in the inclusion and $c(\bx) = 3000\rm{m/s}$ outside, 
see Figure~\ref{fig:Camembert}.
The search space ${\cal C}$ consists of $N = 400$ Gaussian basis functions
centered at the nodes of a $20\times 20$ uniform grid discretizing $\Omega$.
A constant initial guess $c_o(\bx) = 3000\rm{m/s}$ is used.

In the bottom row in Figure~\ref{fig:Camembert} we compare the velocity estimates
obtained with Algorithm~\ref{alg:prowi} (using parameters $L=9$, $q = 4$, $d=n=16$)
with conventional FWI regularized with adaptive Tikhonov regularization similarly to
operator ROM approach, after performing $60$ Gauss-Newton iterations. We observe 
that Algorithm~\ref{alg:prowi} gives a much better estimate of $c(\bx)$ that includes 
a correct reconstruction of both the top and bottom of the inclusion. Conventional FWI 
estimate does not improve much after the $10^{\rm th}$ iteration, indicating that it is 
stuck in a local minimum. Moreover, FWI fails to fill in the inclusion with the 
correct velocity, overestimating it in the upper half of the disk and underestimating it in 
the lower half.

\subsection{Marmousi example}

\begin{figure}[ht!]
\begin{center}
\begin{tabular}{cc}
{\scriptsize True velocity} & {\scriptsize FWI velocity estimate} \\
\hspace{-0.1in}\includegraphics[width=0.463\columnwidth]
{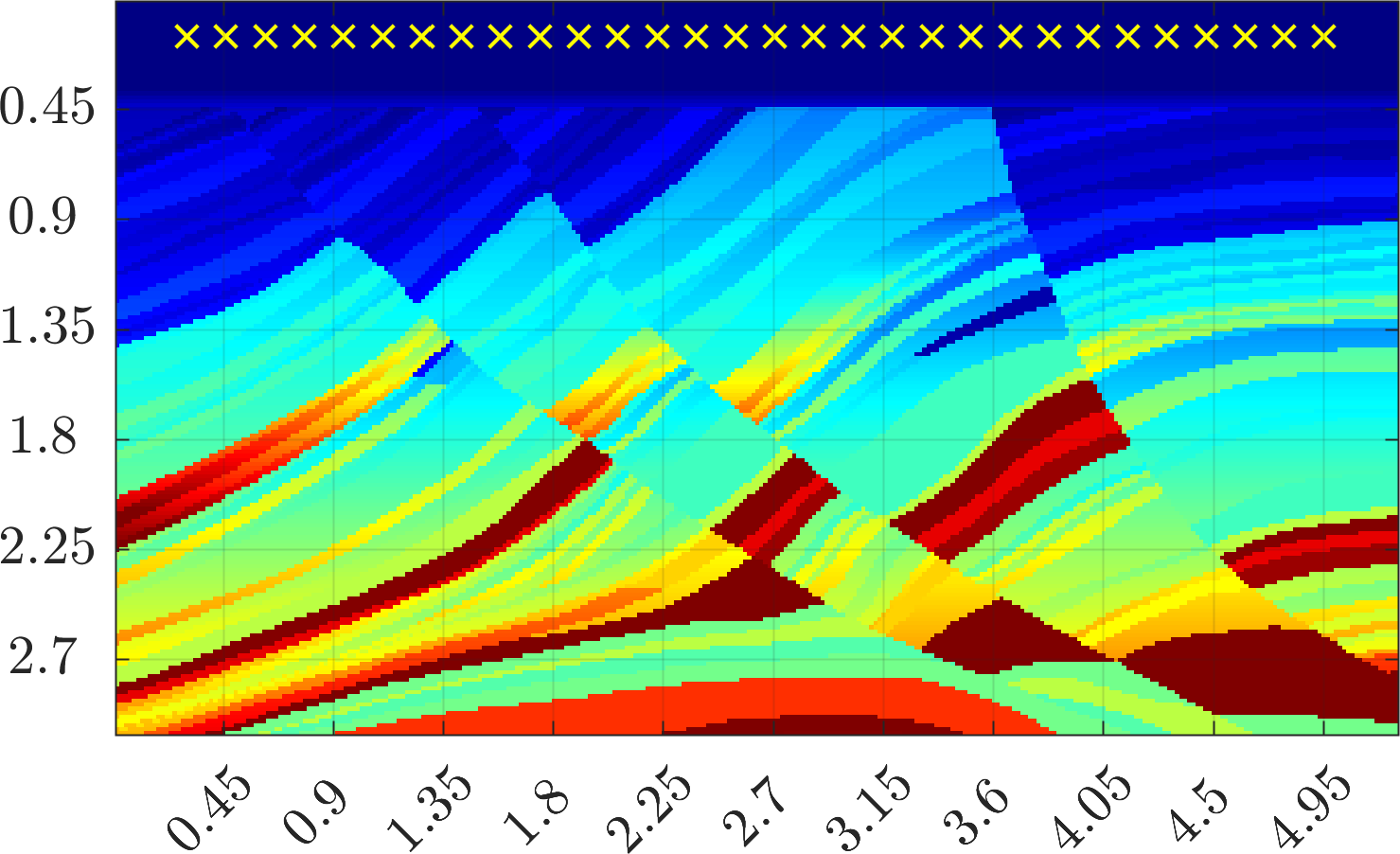} &
\hspace{-0.08in}\includegraphics[width=0.463\columnwidth]
{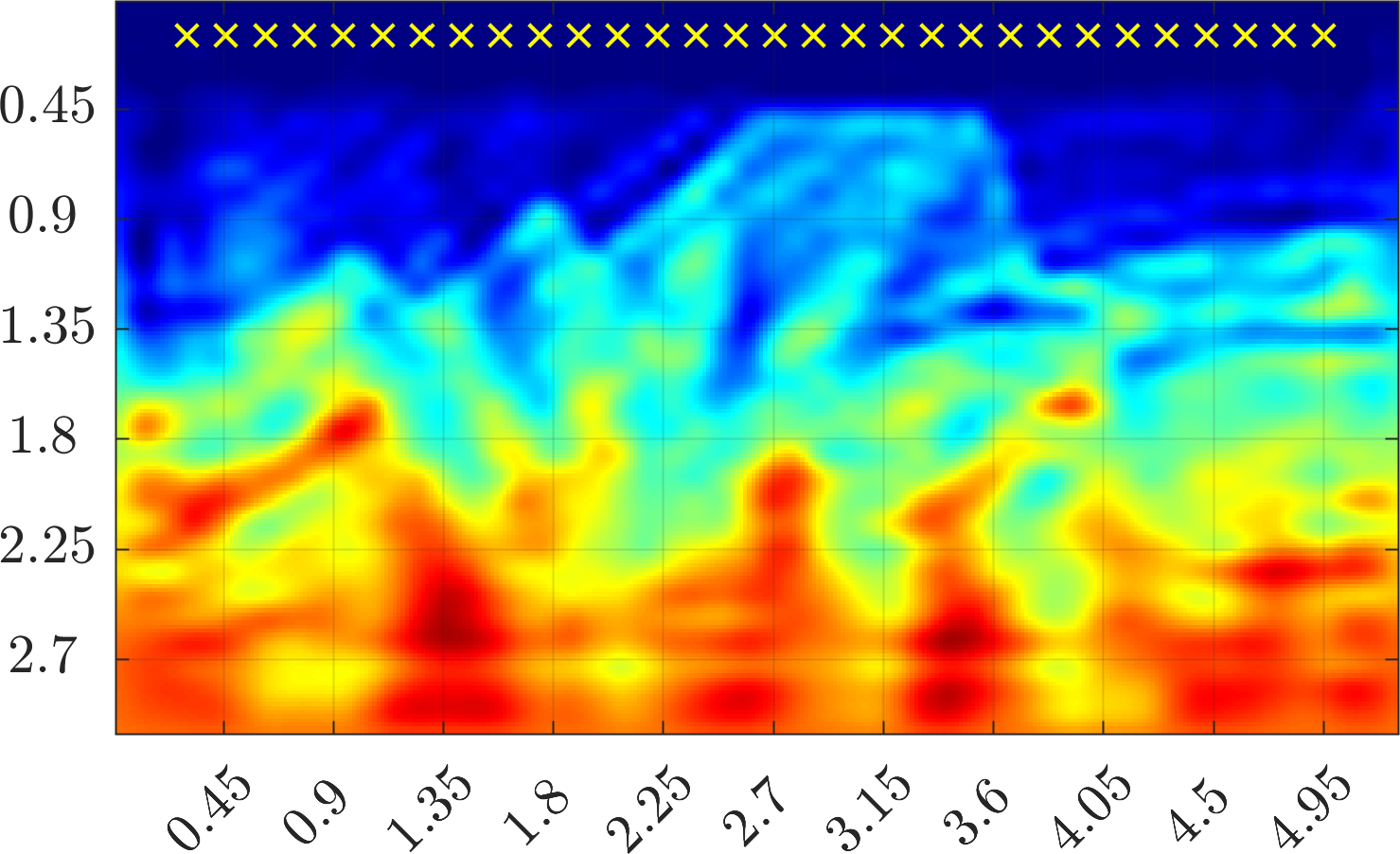}\\
{\scriptsize Refined ROM estimate} & {\scriptsize Operator ROM estimate}\\
 \hspace{-0.1in}\includegraphics[width=0.47\columnwidth]
{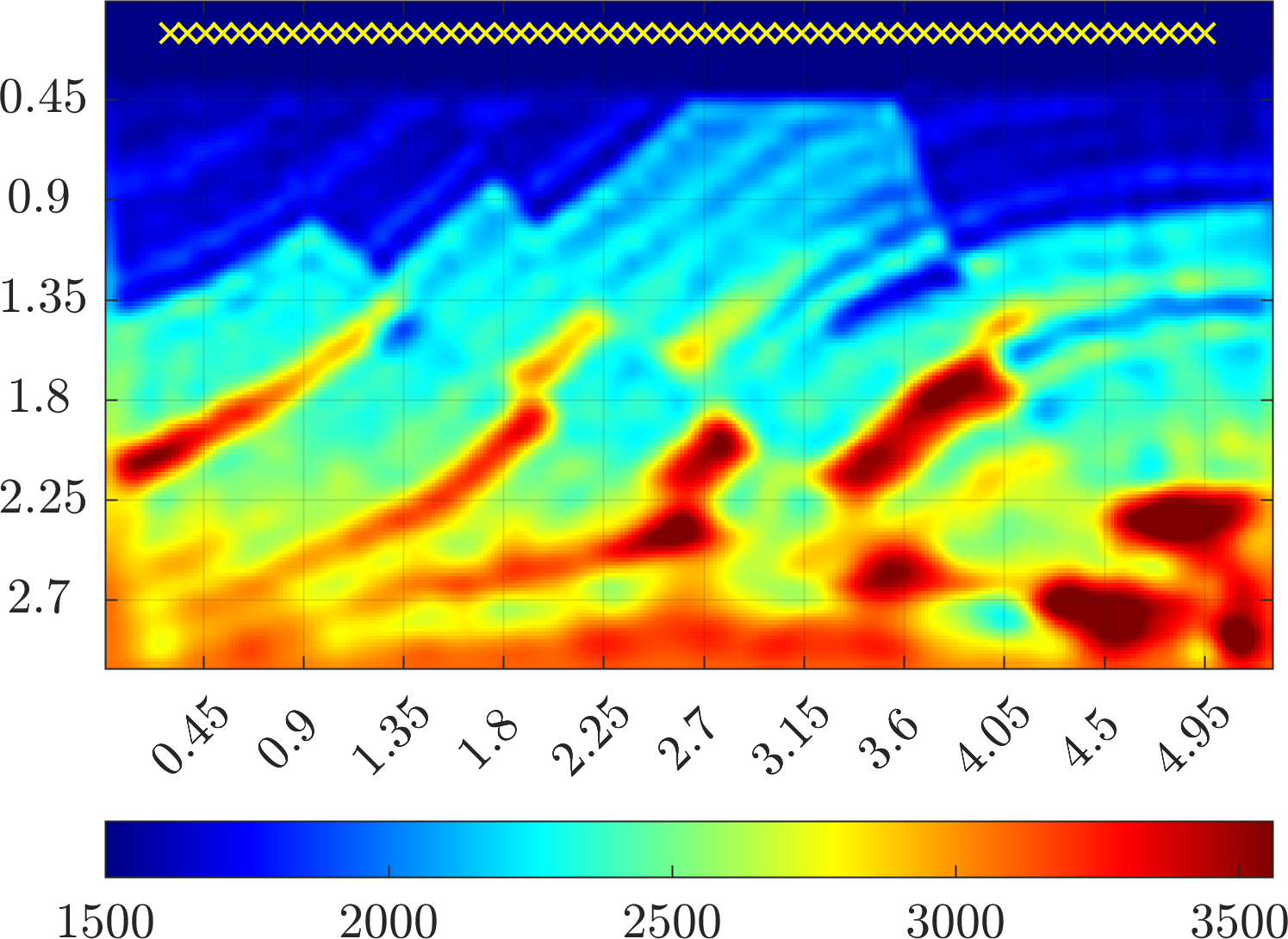} &
\hspace{-0.1in} \includegraphics[width=0.47\columnwidth]
{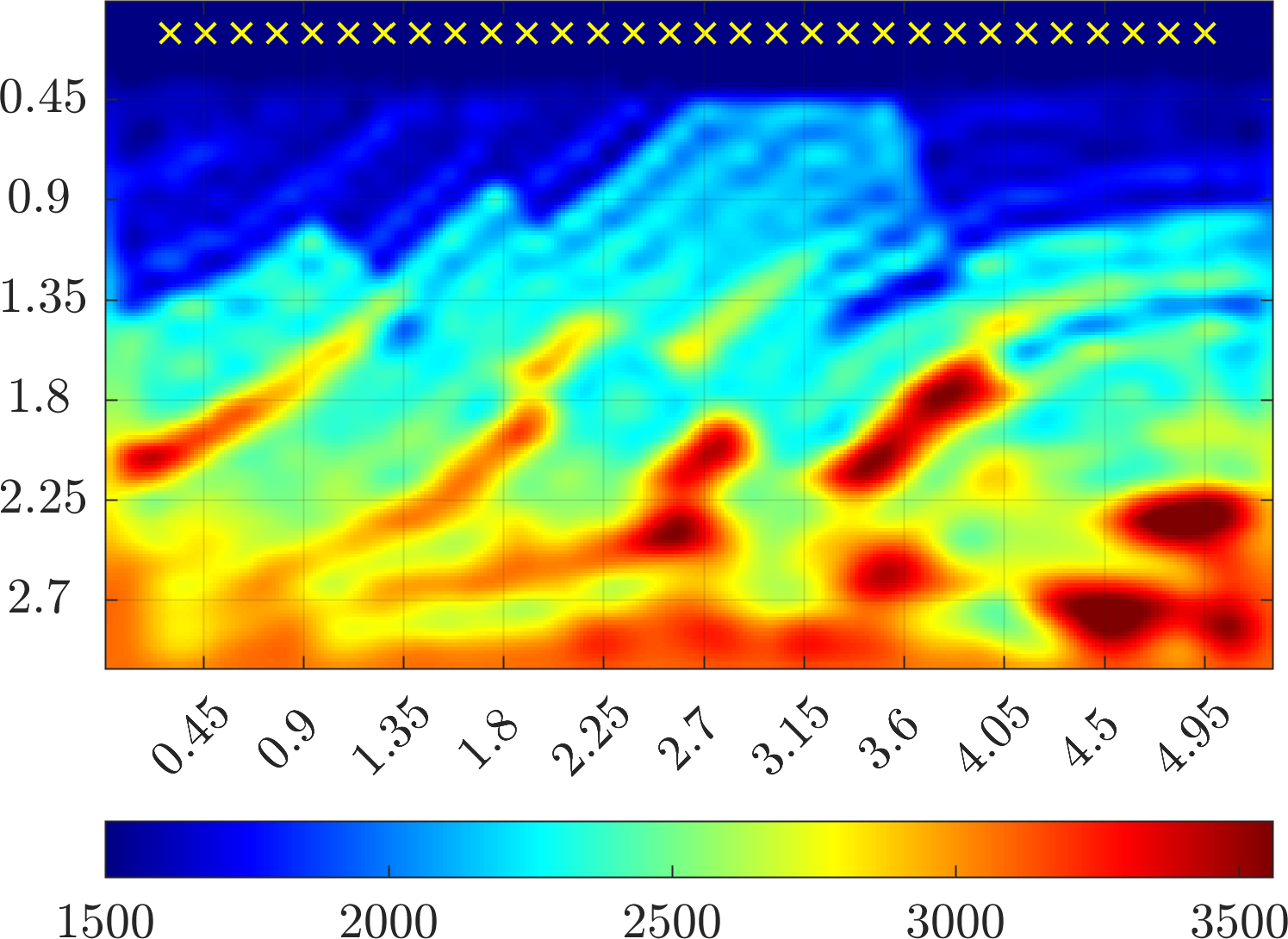}
\end{tabular}
\end{center}
\caption{The section of Marmousi model (top left), velocity estimates 
obtained with Algorithm~\ref{alg:prowi} (bottom row: left $m=60$, right $m=30$)
and conventional FWI (top right, $m=30$). The sensors are shown as yellow $\times$. 
The axes are in $\rm{km}$, velocity in $\rm{m/s}$.}
\label{fig:Marmousi}
\end{figure}

The final example is a section of Marmousi model with water layer down to depth $266{\rm m}$ 
removed. The domain is $\Omega = [0,5.25{\rm km}]\times [0,3{\rm km}]$. 
The search space ${\cal C}$ consists of $N = 1500$ Gaussian basis functions
centered at the nodes of a $50 \times 30$ uniform grid discretizing $\Omega$.
The initial guess $c_o(\bx)$ is a one dimensional gradient in depth.

In the right column of Figure~\ref{fig:Marmousi} we compare the velocity estimates
obtained from data recorded at $m=30$ sensors with Algorithm~\ref{alg:prowi} 
(using parameters $L=6$, $q = 3$, $d=10$, $n=40$) with conventional FWI regularized 
with adaptive Tikhonov regularization after $18$ Gauss-Newton iterations.
We note that the ROM based inversion captures correctly most of the features 
of Marmousi model, while the conventional FWI velocity estimate suffers from a number of
artifacts. We also display in the bottom left plot in Figure~\ref{fig:Marmousi} a refinement 
of operator ROM velocity estimate obtained by injecting more data from $m=60$ sensors
and performing $4$ additional Gauss-Newton iterations of Algorithm~\ref{alg:prowi} 
(using parameters $L=1$, $q=4$, $d=10$, $n=50$) with a refined basis of 
$N = 75 \times 38 = 2850$ Gaussian functions. The resulting refined velocity estimate sharpens 
the boundaries of the features and improves their contrast.

\section{Conclusions}

We presented a novel approach for velocity estimation based on the wave operator ROM.
The ROM is computed from the data and is used to formulate the optimization problem 
for velocity estimation as ROM misfit instead of the conventional FWI data misfit minimization. 
This has a convexification effect on the optimization objective, as shown in a numerical study
of objective topography. As a result, the proposed approach outperforms the conventional 
FWI in synthetic examples such as the ``Camembert'' and Marmousi models.

\section{ACKNOWLEDGMENTS}
This material is based upon research supported in part by the U.S. Office of Naval Research 
under award number N00014-21-1-2370 to Borcea and Mamonov. Borcea, Garnier and Zimmerling 
also acknowledge support from the AFOSR awards FA9550-21-1-0166 and FA9550-22-1-0077.
Zimmerling also acknowledges support from the National Science Foundation under Grant No. 2110265.







\bibliographystyle{seg}  
\bibliography{biblio}

\end{document}